\documentclass[11pt,a4paper]{article}
\usepackage[left=3cm,right=2.5cm,top=2cm,bottom=2.5cm,a4paper]{geometry}

\usepackage[math]{anttor}

\usepackage{amsfonts,amsmath,amssymb,amsthm,calc,mathtools}
\usepackage[utf8]{inputenc}
\usepackage[portuguese]{babel}
\usepackage[T1]{fontenc}

\usepackage{lettrine}

\usepackage{pictexwd}
\usepackage{graphbox,graphicx}
\usepackage{float}

\usepackage{bm}
\usepackage{setspace}

\usepackage{booktabs}

\usepackage{xcolor}
\usepackage[colorlinks]{hyperref}
\hypersetup{
    colorlinks=true,
    linkcolor=blue,
    filecolor=magenta,      
    urlcolor=cyan,
    citecolor=violet,
}

\date{}
\author{Carlos Gomes\\
\small{Escola Secundária de Amarante}\\
\small{\texttt{carlosgomes@esamarante.edu.pt}}\\
}

\title{A Geometria da Regressão Linear}

\begin{document}

\maketitle

\begin{abstract}
\noindent
A regressão linear é um tema normalmente explorado (nas escolas) com recurso a uma calculadora científica gráfica ou software da moda (GeoGebra, por exemplo), ficando os estudantes com a tarefa aborrecida de introduzir números em listas e obter, como recompensa, uma equação que utilizam para fazer previsões num dado contexto. O que aqui se trata é de mostrar o grande valor didático deste problema, mobilizando conhecimentos que os alunos detêm para aclarar, do ponto de vista geométrico, o que está em causa em todo este processo que decorre nos ``bastidores'' da tecnologia.

\end{abstract}

\section{A geometria do problema}


O problema que consiste na determinação da recta que melhor se ajusta a uma dada nuvem de $n$ pontos $(x_i,y_i)$ é tradicionalmente tratado como o problema de encontrar os parâmetros $a$ e $b$ da equação $y=ax+b$ que minimizam a soma
$S = \sum\limits_{i = 1}^n {d_i^2 }$, em que os $d_i$ são as diferenças entre os valores observados e os valores do modelo, isto é,
$
d_i  = y_i  - ax_i  - b
$
(veja-se \cite{gazeta168}).

\begin{figure}[h!] 
 \centering
 \includegraphics[bb=0 0 295 267,keepaspectratio=true, scale=1]{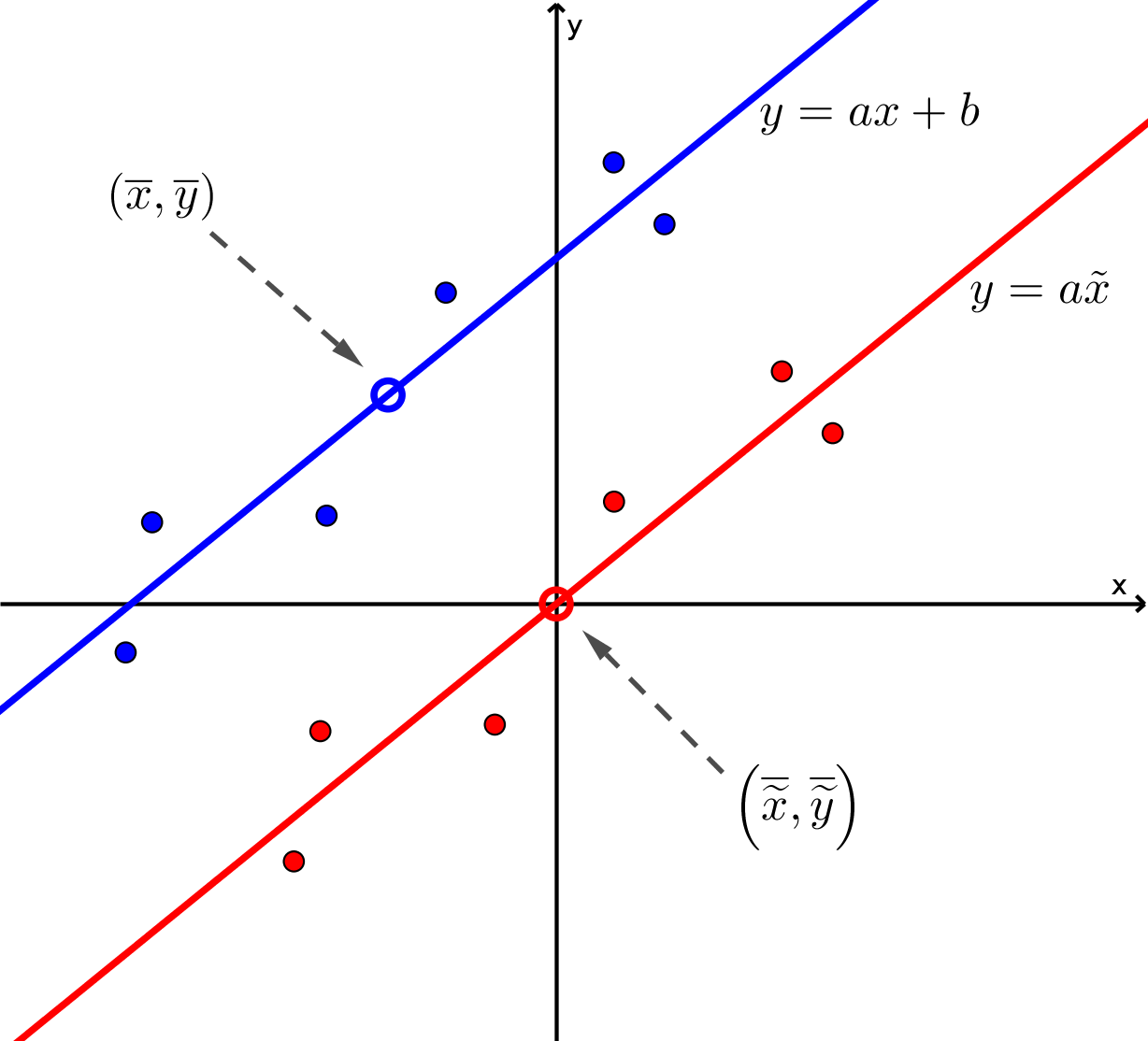}
 \caption{Translação da nuvem de pontos.}
 \label{fig.1}
\end{figure}

Sejam $( {x_1 ,y_1 }),\;({x_2 ,y_2 }),\; \ldots \;,({x_n ,y_n })$ os dados observados (nuvem de pontos na Figura \ref{fig.1}). Para a determinação do parâmetro $a$ (declive da recta), seria ``simpático'' que a nuvem tivesse o seu centro de massa na origem do referencial, isto é, no ponto de coordenadas $( {0 ,0 })$. Isto porque libertar-nos-íamos do parâmetro $b$ da equação da recta, o que parece reduzir a dificuldade do problema, pois, nesta condições, o modelo associado à recta de regressão seria $y=ax$. Para fazer com que o centro de massa da nuvem se desloque para a origem, é suficiente efectuarmos uma translação de toda a nuvem de pontos segundo o vector
$({ - \overline x , - \overline y })$, ou seja, basta subtrairmos o centro de massa $({\overline x ,\overline y })$ a todos os pontos da nuvem. Obtém-se assim uma nova nuvem de pontos da forma $({x_i  - \overline x ,y_i  - \overline y })$ cujo centro de massa é $({0,0})$.

Fazendo $x_i  - \overline x  = \tilde x_i $ e $y_i  - \overline y  = \tilde y_i$, a nuvem sobre a qual o trabalho prossegue será $\left(\tilde x_i ,\tilde y_i\right)$, com $i=1,2,\ldots,n$, cuja recta de regressão tem o mesmo declive que a recta de regressão da nuvem original, em consequência da translação efetuada.

A nova nuvem é constituída por pontos da forma $\left(\tilde x_i,\tilde y_i\right)$ e os pontos da forma $\left(\tilde x_i,a\tilde x_i\right)$, $i=1,2,\ldots,n$, são os pontos sobre a recta $\tilde{y}=a\tilde{x}$, que coincidiriam com os primeiros caso a correlação fosse perfeita. Os $n$ vectores $\vec u_i=\left(\tilde x_i ,a\tilde x_i\right)$ determinados por estes pontos são colineares. Mas aqui, uma mudança de dimensão vai tornar o trabalho mais simples: em vez de considerarmos estes $n$ vectores de dimensão 2, utilizamos os dados organizados em {\bf vectores de dimensão $\bm{n}$}: 
\begin{align*}
{\vec i}&=\left(\tilde x_1,\tilde x_2,\;\ldots\;,\tilde x_n\right),\\
{\vec j}&=\left(a\tilde x_1,a\tilde x_2,\;\ldots\;,a\tilde x_n\right),\\
&e\\
{\vec u}&=\left(\tilde y_1,\tilde y_2,\;\ldots\;,\tilde y_n\right).
\end{align*}
Os vectores ${\vec i}$ e ${\vec j}$ são colineares:
\begin{equation}\label{colinear}
 \begin{split}
	\vec j &=\left(a\tilde x_1,a\tilde x_2,\;\ldots\;,a\tilde x_n\right)\\
 	&= a\left(\tilde x_1,\tilde x_2,\;\ldots\;,\tilde x_n\right)\\
	&= a\,{\vec{i}}\,.
 \end{split}
\end{equation}

Para além do mais, o escalar $a$ em (\ref{colinear}) é precisamente o declive da recta procurada! Assim, determinar $a$ será equivalente a determinar (algo sobre) ${\vec j}$, agora num {\bf espaço de dimensão $\bm{n}$\footnote{Veja o Apêndice no final do artigo para melhor clarificação.}
}.
\begin{figure}[H]
\begin{center}
 \strut\input{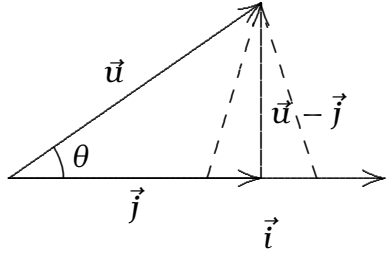}
 \caption{Vetores num espaço de dimensão $n$}
 \label{fig.2}
\end{center}
\end{figure}

Repare-se que
${\vec u}-{\vec j}=\left({\tilde{y_1}-a\tilde{x_1},\;\ldots\;,\tilde{y_n}-a\tilde{x_n}}\right)$
não é mais do que o vector dos resíduos, isto é, o vector cujas componentes são as diferenças entre os dados observados e os dados teóricos da nova nuvem. Ora, o que se pretende é que a norma (ou distância) $\|{{\vec u}-{\vec j}}\|$ seja mínima. Isto só acontecerá se ${\vec u}-{\vec j}$ for normal a ${\vec i}$ (como sugere a figura \ref{fig.2}). Para que tal aconteça, ${\vec j}$ tem de ser a projecção de ${\vec u}$ sobre ${\vec i}$. Logo, o produto escalar de ${\vec u}-{\vec j}$ com ${\vec i}$ tem de ser nulo, retirando-se desta condição o valor do multiplicador $a$, declive da recta de regressão: 
\begin{align}\label{proj}
\nonumber
\left(  {\vec u}-{\vec j}\,\right)  \cdot{\vec i} &= 0\\\nonumber
\Leftrightarrow \left(  {\vec u}-a\,{\vec i}\,\right) \cdot{\vec i} &= 0 \quad \left( \vec j = a\,\vec i,\; de\; (\ref{colinear})\right) \\\nonumber
\Leftrightarrow {\vec u}\cdot{\vec i}-a\,{\vec i}\cdot{\vec i}&= 0 \\
\Leftrightarrow a &= \dfrac{{\vec u}\cdot{\vec i}}{{\|\vec i\|}^{2}}\quad \left({\vec i}\cdot{\vec i} = {\|\vec i\|}^{2}\right)\,.\\\nonumber
\end{align}
Depois de se calcular $a$ através de (\ref{proj}), a determinação do parâmetro $b$ é um simples exercício: dado que $({\bar x,\bar y})$ pertence à recta procurada, ele terá de satisfazer a condição $y=ax+b$. Daqui se retira que $b=\bar y - a\bar x$.

\section{Exemplos de aplicação}

\subsection*{Exemplo 1}

Vejamos a aplicação destes resultados a um exercício típico de um manual escolar.\\

\emph{Existirá alguma relação entre a temperatura e a quantidade de chuva que cai em Amarante? Para responder a esta pergunta vamos comparar num gráfico de dispersão as temperaturas médias (ºC) dos vários meses do ano com a pluviosidade média (mm).}\\

\begin{minipage}[l]{0.5\linewidth}
\centering
\begin{tabular}{@{} c c @{}}
\hline
\emph{Temperatura} & \emph{Pluviosidade}\\\hline\hline
11.3 & 122 \\
12.0 & 108 \\
13.5 & 101 \\
15.2 & 54 \\
17.6 & 44 \\
20.0 & 22 \\
22.2 & 4 \\
22.5 & 6 \\
21.3 & 29\\
18.3 & 80\\
14.2 & 102\\
11.6 & 107\\\hline
\end{tabular}
\end{minipage}%
\begin{minipage}[c]{0.5\linewidth}
\centering
\begin{tabular}{@{} c c @{}}
\hline
\emph{$\sim$Temperatura$\;\backsim$} & \emph{$\sim$Pluviosidade$\;\backsim$}\\
$\vec i$ & $\vec u$\\\hline\hline
-5.3417 &	57.0833\\
-4.6417	& 43.0833\\
-3.1417	& 36.0833\\
-1.4417	& -10.917\\
0.9583	& -20.9167\\
3.3583	& -42.9167\\
5.5583	& -60.9167\\
5.8583	& -58.9167\\
4.6583	& -35.9167\\
1.6583	& 15.08333\\
-2.4417	& 37.08333\\
-5.0417	& 42.08333\\\hline
\end{tabular}
\end{minipage}\\

\vspace*{1ex}

Neste exemplo, a tabela da esquerda é dada e a da direita foi calculada por nós. O centro de massa da nuvem de pontos é $(\bar x,\bar y)=(16.6417,64.9167)$.
Os vectores $\vec u$ e $\vec i$ são as colunas da tabela da direita, depois de efectuada a translação da nuvem original: {\bf são vectores num espaço de dimensão 12}.\\

De acordo com as conclusões da secção anterior, os parâmetros da equação da recta de regressão $y=ax+b$ podem ser calculados do seguinte modo:\\
\begin{minipage}[c]{0.4\linewidth}
\begin{align*}
	a&=\frac{\vec u\cdot\vec i}{{\Vert\vec i\Vert}^2}\\[1,5ex]
	&\approx \frac{-1895.4583}{195.2692}\\[1,5ex]
	&\approx -9.7069\, ,\\[2ex]
	b&= \bar y - a\bar x\\
	&\approx 64.9167+9.7069 \times 16.6417\\
	&\approx 226.4557\,.
\end{align*}
Assim, $y\approx-9.7069x+226.4557$ será a equação da recta de regressão e, com ela, podemos fazer estimativas no contexto do problema.\\
\end{minipage}%
\begin{minipage}[c]{0.6\linewidth}
\flushright
\includegraphics[scale=0.25]{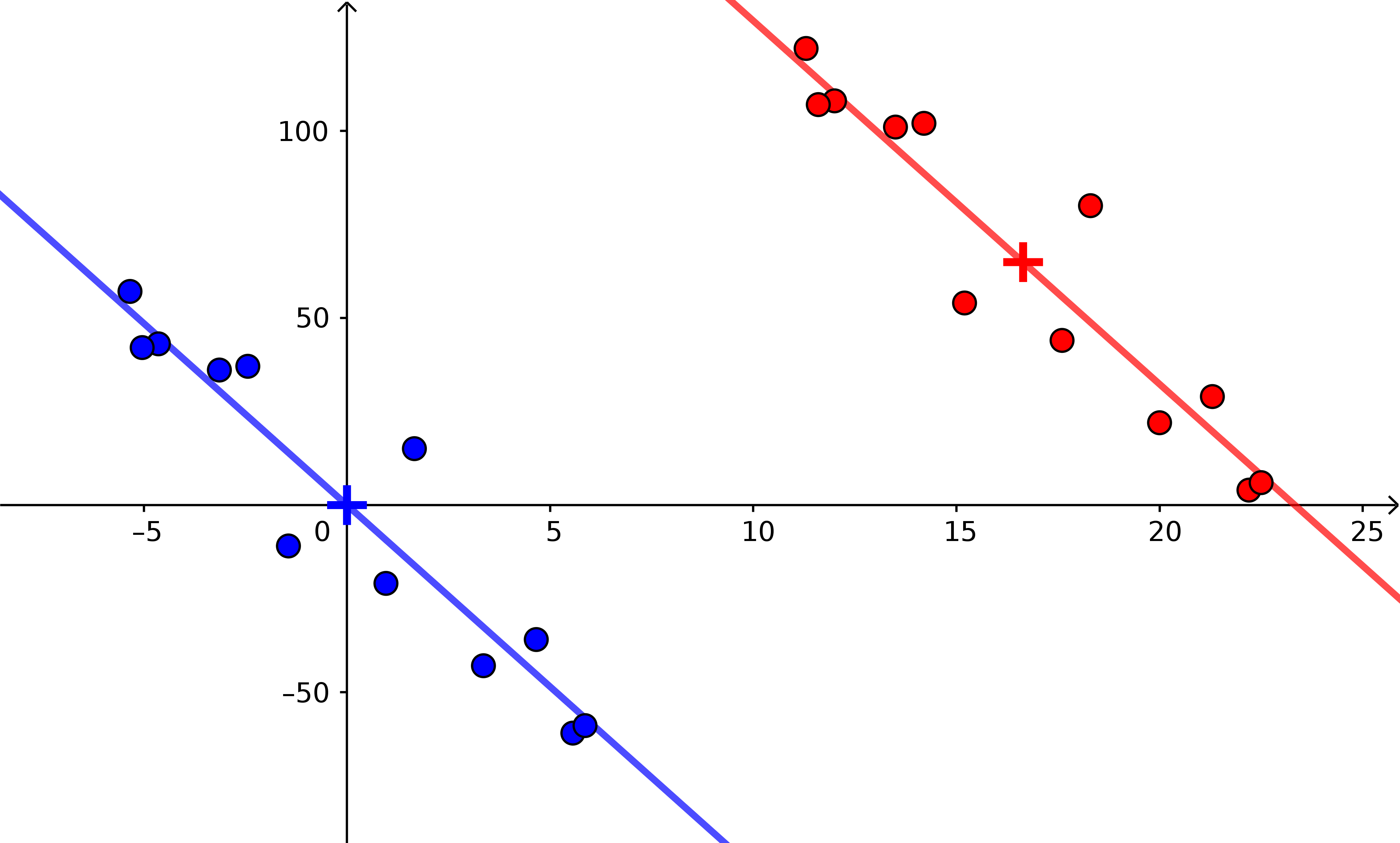}
\label{Fig 3}
\end{minipage}

Note-se que o produto escalar de dois vectores de dimensão $n$ não é mais do que a soma dos produtos das correspondentes componentes desses vectores (uma generalização do que se faz para $n=2$ ou $n=3$ na disciplina de Matemática A no Ensino Secundário), ou seja, se $\vec{a}=\left(a_1,a_2,\cdots,a_n \right)$ e $\vec{b}=\left(b_1,b_2, \cdots, b_n \right)$,

\[
\vec{a}\cdot \vec{b}=a_1\times b_1 + a_2 \times b_2 + \cdots + a_n \times b_n=\sum_{i=1}^{n}{a_i \times b_i}
\;\text{.}\]

Também a norma de um vector de dimensão $n$ é uma generalização da norma de vectores em 2 e 3 dimensões, isto é,

\[
 \lVert \vec{a} \rVert = \sqrt{{a_1}^{2}+{a_2}^{2}+ \cdots + {a_n}^{2}}=\sqrt{\sum_{i=1}^{n}{a_i}^{2}}
\;\;\text{.}\]

Assim, no presente exemplo, $\vec{u} \cdot \vec{i}$ corresponde a efectuar a soma dos produtos dos elementos correspondentes de cada linha da tabela da direita.

\subsection*{Exemplo 2\footnote{Para quem quiser criar uma lição no Geogebra Classroom com este exemplo, seguir para \url{https://www.geogebra.org/m/ncpffvne}}}

Neste exemplo, aplicaremos os conceitos anteriores à construção de um modelo linear do número de infetados pelo novo coronavírus em função do tempo decorrido no período de 8 a 31 de maio. Aqui, o centro de massa é dado pelas coordenadas do ponto $(\bar x,\bar y)=(11.5,29648.583)$ e os vetores $\vec i$ e $\vec u$ {\bf habitam um espaço de dimensão 24} (colunas da tabela da direita).

\vspace*{0.5cm}
\begin{minipage}[l]{0.5\linewidth}
\centering

\begin{tabular}{@{}cc@{}}
\toprule
Nº de dias & Nº de infetados \\ \midrule
67                         & 27268           \\
68                         & 27406           \\
69                        & 27581           \\
70                        & 27679           \\
71                         & 27913           \\
$\cdots$                     & $\cdots$     \\
87                        & 31596           \\
88                        & 31946           \\
89                        & 32203           \\
90                        & 32500           \\ \bottomrule
\end{tabular}
\end{minipage}%
\begin{minipage}[c]{0.5\linewidth}
\centering
\begin{tabular}{@{}cc@{}}
\toprule
$\sim$Nº de dias$\sim$ & $\sim$Nº de infetados$\sim$\\
$\vec i$                   & $\vec u$\\\midrule
-11.5                     & -2380.583       \\
-10.5                     & -2242.583       \\
-9.5                      & -2067.583       \\
-8.5                      & -1969.583       \\
-7.5                      & -1735.583       \\
$\cdots$                     & $\cdots$     \\
8.5                       & 1947.417        \\
9.5                       & 2297.417        \\
10.5                      & 2554.417        \\
11.5                      & 2851.417        \\ \bottomrule
\end{tabular}
\end{minipage}\\[1ex]

O produto escalar é $\vec u\cdot\vec i \simeq 261980$ (soma dos produtos dos elementos de cada linha da tabela da direita). O quadrado da norma do vector $\vec i$ (quadrância de $\vec i$) é ${\Vert\vec i\Vert}^2 = 1150$.

\begin{figure}[H]
\begin{minipage}[c]{\linewidth}
\centering
\includegraphics[scale=0.07]{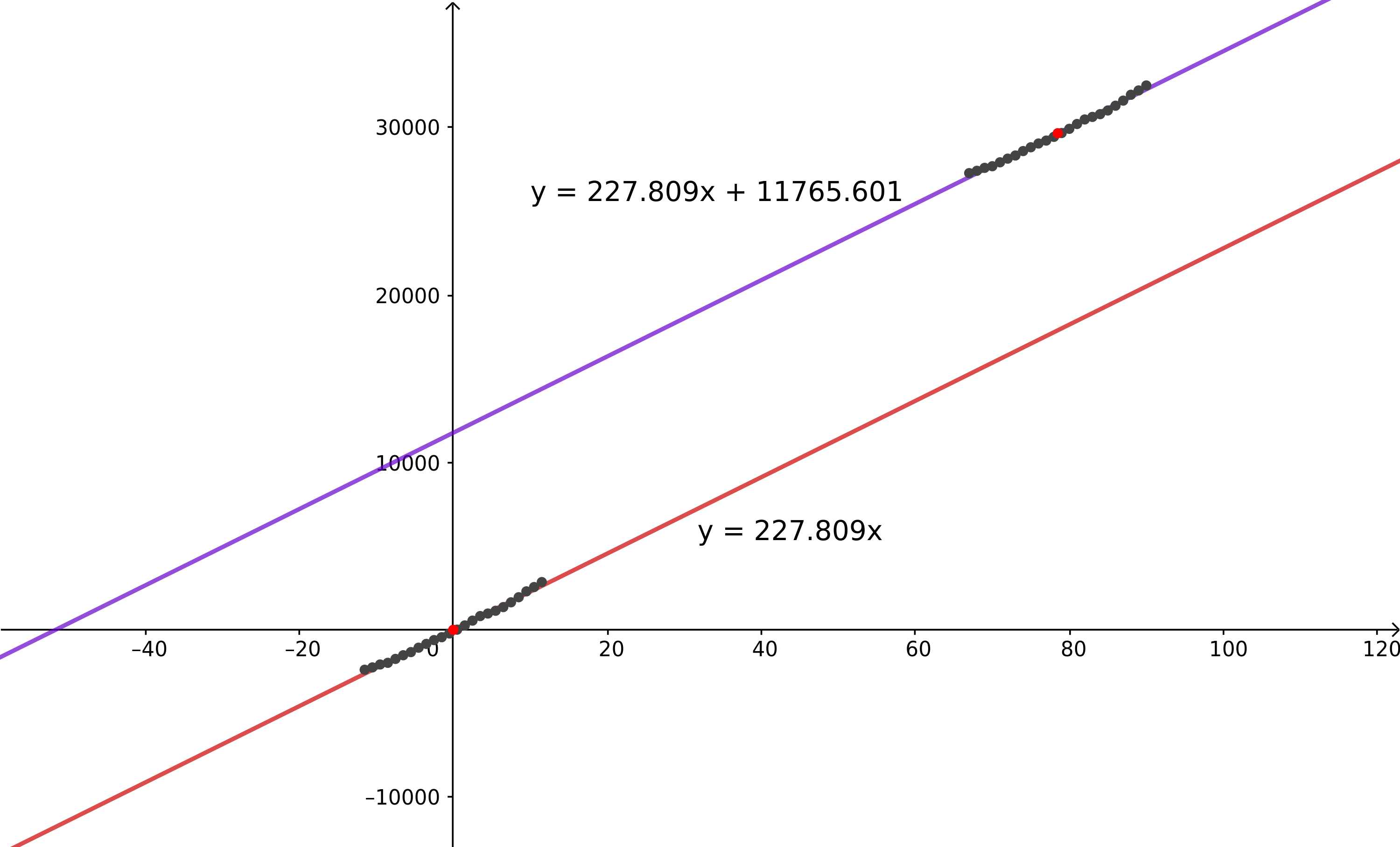}
\end{minipage}
\end{figure}

\vspace*{2ex}

Assim, com $a=\frac{261980}{1150}\simeq227.809$ e $b=\bar y-a\bar x\simeq11765.601$, obtemos a equação da reta mostrada na figura acima.

\section{Coeficiente de correlação linear}

O \emph{coeficiente de correlação} é uma medida que pretende determinar o grau de alinhamento dos dados. Sobre ele costumam ser colocadas duas questões:
\begin{itemize}
\item
\emph{Por que razão varia no intervalo $[-1,1]$?}
\item
\emph{Por que razão a correlação entre as variáveis é tanto mais forte quanto mais próximo de $-1$ ou de $1$ se encontra o coeficiente? Não seria razoável pensarmos que quanto mais próximo de zero mais forte será a correlação, uma vez que ele mede o grau de proximidade dos dados em relação à recta?!}
\end{itemize}

\begin{figure}[H]
\begin{center}
 \strut\input{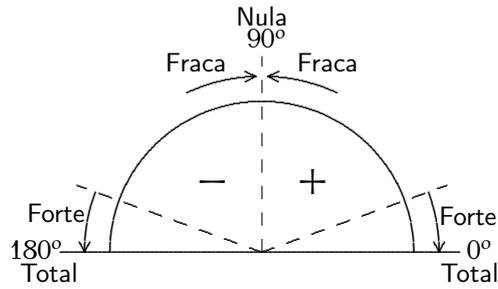}
 \caption{Coeficiente de correlação $\theta$}
 \label{fig.3}
\end{center}
\end{figure}

Repare-se que o coeficiente de correlação, sendo uma medida do alinhamento dos dados, deve estar relacionado com o ``grau de colinearidade'' entre os vectores $\vec u$ e $\vec i$, referentes aos dados transladados \footnote{\scriptsize{A correlação não depende da nuvem que se considera, uma vez que a operação de translação efectuada à nuvem inicial garante a manutenção das relações entre os dados observados e os teóricos.}}. E uma forma natural de medir este ``grau de colinearidade'' é estudar o ângulo $\theta$ que $\vec u$ e $\vec i$ formam entre si (ver figura~\ref{fig.2}).\footnote{\scriptsize{Em tudo o que se segue pode-se substituir a unidade \emph{grau} por \emph{rad}}}
Assim, $\theta$ poderia ser usado com legitimidade como medida do grau de alinhamento dos dados, ou seja, como coeficiente de correlação. O diagrama da figura \ref{fig.3} resume a variação deste coeficiente de correlação.

Visto que
$\displaystyle{\cos \theta = \frac{\vec u \cdot \vec i}{{\|\vec u\|}\,{\|\vec i\|}}\;,}$
$\theta$ pode ser obtido através de
\begin{equation}\label{eq:teta}
\theta = \arccos \Bigg(\frac{\vec u \cdot \vec i}{{\|\vec u\|}\,{\|\vec i\|} }\Bigg).
\end{equation}
No exemplo 1 da secção anterior, o coeficiente de correlação $\theta$ é 
\[
\theta = \arccos \left( \frac{\vec u \cdot \vec i}{{\|\vec u\|}\,{\|\vec i\|} }\right) =
\arccos \left( \frac{-1895.4583}{143.7391 \times 13.9739 }\right) =160.68^o\; (forte\;Negativa?).
\]
e no segundo exemplo, $\theta=\arccos\left(\dfrac{261980}{262579.265}\right)=\arccos(0.998)\simeq 3.62º$ (Muito forte, positiva?).\\

No entanto, na literatura sobre o assunto, $\theta$ é convenientemente substituído pelo seu cosseno (porquê?), e assim se compreende a sua variação tal como encontramos nos manuais:

\[
0^o\le{}\theta\le{}180^o \Rightarrow -1\le{}\cos\theta\le{}1 \Leftrightarrow -1\le{}\displaystyle{\small\frac{\vec u \cdot \vec i}{{\|\vec u\|}\,{\|\vec i\|}}} \le{}1 .
\]

\vspace*{2ex}

Uma fórmula que normalmente acompanha os manuais para determinar o valor do coeficiente de correlação, $r$, é

\begin{equation}\label{eq:r}
\displaystyle{r = \frac{{\sum_{i=1}^{n} {x_iy_i - \frac{{\left( {\sum_{i=1}^{n} x_i } \right)\left( {\sum_{i=1}^{n} y_i } \right)}}{n}} }}{{\sqrt {\left( {\sum_{i=1}^{n} {x_i^2 }  - \frac{{\left( {\sum_{i=1}^{n} x_i } \right)^2 }}{n}} \right)\left( {\sum_{i=1}^{n} {y_i^2 }  - \frac{{\left( {\sum_{i=1}^{n} y_i } \right)^2 }}{n}} \right)}}}}\,.
\end{equation}

Sendo (\ref{eq:r}) equivalente a 

\begin{equation*}
  r =
  \frac{ \sum_{i=1}^{n}(x_i-\bar{x})(y_i-\bar{y}) }{%
        \sqrt{\sum_{i=1}^{n}(x_i-\bar{x})^2}\sqrt{\sum_{i=1}^{n}(y_i-\bar{y})^2}}\,,
\end{equation*}

fica estabelecida a igualdade

\begin{equation*}
\displaystyle{r=\frac{\vec u \cdot \vec i}{{\|\vec u\|}\,{\|\vec i\|}}}=\cos\theta\,.
\end{equation*}

\section{Apêndice}

A interpretação geométrica que se explora neste texto tem como elemento essencial a translação da nuvem de pontos original para uma nuvem de pontos com centro de massa na origem do referencial. Esta operação faz com que as seguintes condições se verifiquem\\
\[
 \sum_{i=1}^{n}{\widetilde{x}_i}=0 \quad \text{e} \quad \sum_{i=1}^{n}{\widetilde{y}_i}=0
\;\text{.}\]

Reescrevendo-as, ficamos com

\begin{align*}
 \sum_{i=1}^{n}{\widetilde{x}_i}=0 & \Leftrightarrow 1 \times \widetilde{x}_1 + 1 \times \widetilde{x}_2 + \cdots + 1 \times \widetilde{x}_n = 0 \\
    & \Leftrightarrow \boxed{\vec{w} \cdot \vec{i}=0}\\
    & {\centering\text{e}} \\
 \sum_{i=1}^{n}{\widetilde{y}_i}=0 & \Leftrightarrow 1 \times \widetilde{y}_1 + 1 \times \widetilde{y}_2 + \cdots + 1 \times \widetilde{y}_n = 0 \\
    & \Leftrightarrow \boxed{\vec{w} \cdot \vec{u}=0}\;\text{,}
\end{align*}

que, do ponto de vista geométrico, permitem afirmar que os vectores $\vec{i}$ e $\vec{u}$ (e, consequentemente, $\vec{j}$) são perpendiculares ao vector unitário $\vec{w}=\left(1,1, \cdots, 1\right)$. Assim, $\vec{i}$, $\vec{j}$ e $\vec{u}$ habitam o hiperplano de dimensão $n-1$, normal ao vector unitário $\vec{w}$. Este facto não altera a argumentação seguida pois no hiperplano de dimensão $n-1$ continuamos a querer reduzir ao mínimo a norma de $\vec{u}-\vec{j}$ e a condição continua a ser a ortogonalidade deste vector a $\vec{j}$. No caso em que a amostra observada é constituída apenas por dois pontos, $\vec{i}$, $\vec{j}$ e $\vec{u}$ são colineares e a correlação é perfeita, como seria de esperar\footnote{Para a situação em que $n=3$, pode manipular e descarregar a animação GeoGebra em\\ \url{https://www.geogebra.org/m/muxygsbz}}.

\section{Conclusão}

Ao longo dos anos, o tema da regressão linear tem sido tratado nas nossas escolas, quase exclusivamente, como uma manipulação de fórmulas, à qual a tecnologia veio retirar algum desse desprazer salvando, por um lado, os alunos dos cálculos fastidiosos, mas atirando-os, por outro, para uma cegueira determinada pela calculadora gráfica. O que aqui se quis mostrar foi que essas abordagens tradicionais ao tema podem, com enormes vantagens, serem substituídas por uma abordagem geométrica sólida, coerente e palpável, em que a única novidade (mas não surpresa) reside na generalização de conceitos de geometria analítica a espaços de dimensão superior a três. Para além disso, abre também espaço à compreensão dos ``bastidores'' da calculadora gráfica, permitindo que os alunos olhem para ela como uma biblioteca de algoritmos que podem compreender e até criar.

\end{document}